\documentclass{paper}[12pt] 

\usepackage{amssymb}   
\usepackage{amsthm}
\usepackage{amsmath}
\usepackage{amsfonts}
\usepackage{latexsym} 
\usepackage{eucal}
\usepackage{indentfirst} 
\usepackage{graphicx} 
\usepackage{verbatim} 

\begin{document}

\overfullrule=0pt
\baselineskip=24pt
\font\tfont= cmbx10 scaled \magstep3
\font\sfont= cmbx10 scaled \magstep2
\font\afont= cmcsc10 scaled \magstep2
\title{\tfont Using Mathematical Maturity to Shape our Teaching, our Careers and our Departments }
\bigskip
\bigskip
\author{Thomas Garrity\\  
Department of Mathematics and Statistics\\ Williams \\ Williamstown, MA  01267\\ 
email:tgarrity@williams.edu \thanks{This a written version of a talk originally given at the conference  ``The Art of Teaching Mathematics" at Harvey Mudd in June 2007.  Other versions were given at an MAA Northeastern meeting at Framingham in November 2007, at the Park City Mathematics Institute (PCMI) in June 2008, at a Project NEXT workshop at an MAA Northeastern meeting at Bradley University in November 2008, and at  MAA dinners at BYU in March 2008 and at Norwich University in April 2008.  To all who heard and commented, I owe thanks.  I would also like to thank  Jon Jacobsen and Michael Orrison for the invitation to the Harvey Mudd conference and to John McCleary, Frank Morgan and Lori Pedersen for many useful comments. An abridged version is in \cite{Garrity11}.  After this was written, there appeared \cite{Krantz11},}}

\date{}

\maketitle

From 2006 to 2013, I was the director of  the Williams College Project for Effective Teaching (Project PET, which is now called First Three).  While PET  does a number of things, its main goal is to help mentor new faculty in all disciplines, hopefully helping  them  to thrive at Williams as both teachers and as scholars.  At a college that prides itself on teaching, somehow I was the teaching guy.  It was a great job, since I  spent a lot of time talking to ambitious smart young people; I've also learned an amazing amount about other disciplines and departments.  One thing that has truly surprised me is that no other discipline has a concept analogous to our mathematical maturity (with the possible exception of departments teaching foreign languages, where there is the notion of ``language readiness')\footnote{It is possible that the  term ``threshold concept, '' recently coined by educators Mayer and Land, is attempting to describe analogs of mathematical maturity to other contexts.  Certainly, though, mathematicians have been using the term ``mathematical maturity'' for decades.}.  I am sure that a professor of history, say, is more Ôhistorically matureÓ than a first year student, but it is not the rhetoric that they use.  In fact, it takes some effort to   even  give people an idea of what we mean by mathematical maturity.   But it is a vital notion.  

I have  come to believe that we can use ideas of mathematical maturity to shape our teaching, our careers and our departments.  In fact,  I suspect that many of the problems in teaching stem from the professor (or teacher) teaching at one level of intellectual maturity while the students are at a different (lower) level.  For example, the difficulty in teaching beginning calculus often stems from  students not being fluent in high school algebra. I'm guessing that there are analogs in other disciplines. You could imagine a philosophy professor talking about the details of the synthetic a priori  in Kant  to students who have never questioned basic empiricism.

Using the rhetoric of mathematical maturity might also  help shift  the standard arguments of reform versus non-reform teaching methods to an  emphasis  on how we should all be trying to foster our own levels of mathematical maturity, from "the cradle to the grave".  For example, I do not think of math quite the same now as I did 20 years ago, even though 20 years ago I was already on the faculty at Williams.  I hope I will have even more of an understanding of math in 20 more years, when I am possibly thinking about retirement.

First, though, what exactly is ``mathematical maturity"?  Most of us probably think of this as primarily being about understanding proofs.  We see  huge leaps in understanding in our students when they finally get it.  At Williams, this most often happens in the first course in real analysis or in the first course in abstract algebra.  A student is struggling to just pass, and then is suddenly  doing A work.  

Thus, while I can't quite define mathematical maturity, I've  seen it ``in action.''   Many of us probably remember our own leap into this type of mathematical maturity.  For me, it happened in the first month of college.  Bruce Palka was setting up a new honors program in math at  the University of Texas.  Our text was Spivak's {\it Calculus}, and the course was  pitched at the level of this book. I found the course  hard and had at first no real clue about what was going on.  How I got through the first homework problem set, I still don't know.  On the second, though not deliberately, I basically copied off of Michael Lacey (who was far more mathematically mature than I;  in fact, I suspect he still is).  He didn't bother to show up for a meeting for the third homework assignment (I don't blame him),   leaving me with my equally ignorant colleagues.  Then it happened.  It was a Thursday,   David Bowie's Diamond Dogs was echoing in my head, and I got it.  Suddenly problems that seemed impenetrable  became intriguing.  I could work the problems.  This was one  of the greatest days of my life.

But there are kinds of mathematical maturity that come even before doing proofs.   Certainly one of the big stumbling blocks for many beginning college students is their lack of fluency in high school algebra.  Those who speak this language can get through almost all of the first few years of college math.  Those who don't, struggle with calculus.  Speaking this language is also a  type of mathematical maturity.  Maybe most of you found your first exposure to basic algebra easy.  I didn't.  We had an ambitious seventh grade teacher who in the middle of the year tried to introduce some algebra to us.  I suspect it was just solving a single linear equation in one variable.  I just didn't get it.  All I remember is that there was an equation, with an equal sign, then on the next line a new equation, still with an equal sign.  I had no clue as to what was going on.  We soon returned to the normal stuff,  coming back to this algebra at the end of the year, at which time it all made perfect sense.  I had become high school algebraically mature.   Note that I experienced no moment of epiphany.

Can this story of mathematical growth  be extended and fleshed out, so that all of us can  foster our own personal mathematical maturities, tracing our own mathematical paths from nursery school to the retirement home? At each stage there are both mathematical facts to learn and levels of mathematical maturity to attain.  For example, by the end of elementary school, students should be comfortable with basic arithmetic (including fractions), know how to represent quantitative data (charts and graphs), have a sense of magnitude (e.g., know that 38 times 43 is about 1600) and be able to recognize basic geometric objects.  The level of mathematical maturity is for students to recognize that the above are not just random facts and techniques but part of larger logical structure.  Sixth graders need to know not proofs but should expect that  math should make sense. Almost all elementary school teachers know the mechanics;  fewer know  the underlying reasonings.

By the end of secondary school, students  should know basic algebra  (able to comfortably manipulate algebraic equations),  should know basic functions such as trig functions, exponentials and logarithms, should know the beginnings of Euclidean geometry (primarily as an example of the axiomatic method) and should know some basic counting and probability theory.  For  mathematical maturity, they should recognize, as in elementary school, that  math is not a bunch of facts but a logical whole built on proof (from  trig identities to  Euclidean geometry).  

By the end of the first two years of college, our budding mathematicians should know basic calculus (at least through multi-variable calculus) and linear algebra (including how to solve big systems of linear equations via matrices and, even more importantly, the need for abstract vector spaces).   In calculus, mathematical maturity means an understanding of  the practical uses of 		     
calculus, the idea of the  derivative and the integral, and an intuitive understanding of the Fundamental Theorem of Calculus.  For linear algebra, it means seeing   when 
 a problem can be reduced to linear algebra, coupled with the recognition that the problem can then (often) be solved.

For the last two years of college, people should know the beginnings of  real analysis, abstract algebra and possibly point set topology, complex analysis,  or differential geometry.  For mathematical maturity,  our college senior math major should  know how to recognize and produce fully rigorous proofs.  (This is 	  	     far more important than any knowledge picked up in classes.)

In the first few years of graduate school, our young mathematician should be learning the broad outline of the first 50 years of 20th century mathematics, including measure theory, homology and cohomology theory, abstract algebra 
	and complex analysis.  More important, and more relevant for mathematical maturity, is to develop the  ability to learn mathematics quickly.

	For the last few years of graduate school, the goal  is to produce a thesis,  to  become an expert in a narrow branch of mathematics and along the way to learn how to deal with the frustration of burrowing into a narrow problem.

In the first  few years after the  Ph.D., scholars should be examining the mathematics near their thesis areas.  In terms of mathematical maturity, these new faculty should be increasingly developing their own personal view of mathematics. Note that I am now talking about young  math professors and thus people who are extremely mathematically mature as compared to the general population.

This sketch of the different levels or types of mathematical maturity can be continued to those of us who are in mid-career, near retirement and even to those of us in the nursing homes.  Who among us would not like to end our days, at 102, in a hospital bed, surrounded by a loving spouse, children, and grandchildren (maybe with one great great grandchild, an infant in arms) and have as our last thought our last moment of insight into mathematics?

How can we use the idea of mathematical maturity to influence our teaching, our careers and our departments?  
For teaching, this is fairly clear.  We can and should use the idea of mathematical maturity to help shape  our classes.  Too often our classes are reduced to the mere teaching of technique (which most students think of as all of mathematics).  We should always bring in the big picture and  keep emphasizing, over and over, what is the goal.   For me,  in each lecture there should be a clearly stated punchline, which in turn should be linked to the goals of the unit, the goals of the semester, and to mathematics overall, all at the age appropriate level of mathematical maturity.

This will not only prevent the  mindless recitation of facts and techniques but also  help  to determine  what is important and what is mere detail.  For example, in calculus, we can see that the Fundamental Theorem of Calculus is more central than  integration by parts, which in turn is more central than integrating inverse trig functions.

This  is also where research can influence our teaching, which rarely happens.  Back in the mid 1980s,    Ken Hoffman  became  the first math lobbyist, after a successful career as a mathematician at MIT.  Back then I heard him speak once  about trying to get an increase in government math  funding.   He said that all congressmen and for that matter everyone else in government thought that math was important.  On that, none had any doubts.  What they didn't know was that math was still going on. Those of us in the audience had at that moment a certain smug self-satistfication, looking down on the philistines in Washington.  But then Hoffman asked who was to blame, pointing an accusing  finger at the audience.   After all,  how many students are in math classes every day?  How could they possibly  know that math is still going on if not told by their teachers, their professors?   Yes, we are to blame.  Mention research.  Not in abstract terms but as how it enters into the classroom. 

For example, in the late 80's as a young postdoc, I cared very much about algorithms for factoring multivariable polynomials over the complex numbers.  This led me in my classes to emphasize factoring as a key type of  problem.  Did I go on and on, with all the details of then current algorithms?  Of course not.  I just mentioned, whenever factoring came up, that if you bump up the number of dimensions then you would be doing current work that was important.  This can be done no matter what type of research you are doing.  Simply occasionally mention how the topic in your given class can be slightly generalized to current mathematics.  

Mathematical maturity should also be used to help guide each of us in our careers.  Most readers of the Notices have PhDs.  A large percentage of the readers work at schools with fairly heavy teaching loads.  These are folks who look sheepish (almost ashamed) if you ask them about their research, despite the fact that if they are teaching four or more courses a semester, it almost certainly precludes time for research.   There has to be another way for these talented people to continue to foster their mathematical maturity and not get stuck explaining pre-calculus over and over again for the next fifty years. Here is one experiment that a group of us are currently trying.  In the summer of 2008, at the IAS/Park City Mathematics Institute (PCMI), I ran the Undergraduate Faculty Program.  The goal was to show faculty from primarily undergraduate institutions how to teach a course in algebraic geometry.   Thus I was lecturing to adult mathematicians.  Everyone knew that   one has to work problems  to understand a new area of math.  But if they were like me, they would look at a given problem, say to themselves that they could work it (sadly without really putting pencil to paper)  and after a few days be overall lost.  At the same time, I realized that my biggest weakness as a teacher is in coming up with good problems. So, instead of assigning problems, we decided to come up with our own.  Our initial goal was that by the end of the three-week institute, we would have a nice collection of problems.  It rapidly became clear, though, that we could turn these problems into a book of problems, an introductory text for algebraic geometry, which has now appeared  \cite{Garrity 13}.  What is more important for this article is that we now have an informal network of collaborators.  The people involved teach at schools with widely varying  teaching loads.  What we had in common was serious interest in teaching.  A possible model would be to form such informal networks among young mathematicians with the goal of collaborating  on possibly expository  work. (This is also a reason for people who care about teaching and research to attend  the Undergraduate Faculty Program at next summer's PCMI (https://pcmi.ias.edu).)

What about our departments?  We all know of departments where people stay in their offices, leaving only  to  teach their classes.  Such departments are only as strong as each individual member and are particularly toxic to junior folks.  Certainly there are some straightforward methods, such as a weekly department colloquium or daily lunch crowd.  For places with heavy teaching loads, I would suggest that the weekly department colloquium be mainly local people talking, as opposed to  outside speakers.   These talks certainly do not have to be cutting edge research.  After all, how many of us know the cutting edge research in an area of math from 1980, say, that is far from our own research.

There are more radical ideas.  At Williams we have experimented with something we call DADA math and with mathematical orgies (which are variations on what happened at Brown a few times back in the eighties).  For DADA,  the motivation is indeed from the artistic movement DADA.  The idea is to get  people in the department talking math with each other.  We start by sitting around a table.  In a ceramic bowl are slips of paper, each with a classification number from Math Reviews.  We pull three slips out, at random, and write down on a blackboard the corresponding topic.  For example, if we pulled out 17D, 51M and 55R, we would write down the topics ``Other nonassociative rings and algebras'', ``Real and Complex Geometry'' and ``Fiber Spaces and Bundles.''  We then start free associating, trying to find an interesting question linking all three topics.  After 15 to 20 minutes, we do it again.  The dream would be to have a true research paper come out of this process.  That has not yet happened.  I think that the most successful DADA was the first, when we actually came up with a clean conjecture that we later discovered had been made by Mahler back in the 60s.  Still, we have had a good time and have learned new mathematics.  For example, at one we came across the term ``reverse mathematics'', something none of us had heard of.  Susan Loepp speculated that this was doing mathematics walking backward while emitting a beeping sound, for safety's sake.  What actually happened is that almost everyone, within a day or two, googled reverse mathematics and got to learn a bit about a fascinating part of mathematical logic.  Overall, our experiments with DADA math have met with mixed success (perhaps partly because at more recent DADA events, beer was present,  resulting in more humor but less insight.)

Now for the orgies. In the 1980s, Brown had a few all day sessions where many of the faculty, graduate students and a few outside experts spent a day at a series of talks, called {\it orgies}, learning about some area of algebraic geometry.   At Williams, starting in 2002, we have altered this tradition a bit. We have a January term for experimental classes, giving the faculty some breathing room, as only half of us  teach during this term on any given year.  Roughly each year we choose a topic and spend a day or afternoon lecturing to each other on it.  Here we do not bring in outside experts.  No one is assumed to be an expert.  Our first one was called the ``Enumerative Orgy.''  The goal was to discuss the links between string theory and enumerative problems in algebraic geometry.  The second, in January 2006, was the ``Semeredi Festival'', where we went through a survey article by Terry Tao \cite{Tao05}.  The third, in 2007, was ``Poincairepolloza'', where we discussed the proof of the Poincaire conjecture.  In 2008, we had``Stochastic Fantastic Day'', where we used David Mumford's ``Dawning of the Age of Stochasticity'' \cite{Mumford00} as a springboard to talk about the nature of randomness.  Most recently, in 2010, we had ``Transferrific Day'', using David Ruelle's ``Dynamical Zeta Functions and Transfer Operators'' \cite{Ruelle02} as a template for the day.  These are great ways to get people in the department to talk to each other about mathematics (instead of administrative issues).  Here is the best story so far from these orgies.  It happened during the Semeredi Festival.  Olga (Ollie) Beaver (who recently passed away)  was one of the speakers.  You have to know that Ollie was  the type of mathematician who was close to a perfectionist.  Making a mistake in a lecture almost caused her physical pain.  On the day of the orgy, Ollie got up and gave a fine thirty minute lecture.  Only later did we hear the real story.  The type of math that she was to speak on was far from her area of expertise.  In the weeks before the orgy, she struggled to understand her assigned section.  She just didn't get it.  She spent nights without sleep.    On the morning of the orgy, she had nothing.  Then, during the talk right before hers, as she sat in the back of the room, trying to remain calm, she had a blinding flash of insight.  She got it, and ended up  giving a polished lecture.  What is important is that by putting herself on the line, she gave herself the opportunity of have this insight.  Though the most senior member of our department, she still took the chance to  become more mathematically mature.  

Hopefully we can all be so lucky.


\begin{thebibliography}{99}



\bibitem{Garrity11}  T. Garrity, Using Mathematical Maturity to Shape Our Teaching, Our Careers, and Our Departments, {\it Notices of the American Mathematical Society}, volume 58, issue 11 (2011), pp. 1592-1593.

\bibitem{Garrity 13}  T. Garrity, R. Belshoff, L. Boos, R. Brown, J. Douilhet, C. Lienert, D. Murphy,  J. Navarra-Madsen, P. Poitevin, S. Robinson, B. Synder, C. Werner, {\it Algebraic Geometry: A Problem Solving Approach}  American Mathematical Society, Student Mathematical Library, Vol. 66, 2013.



\bibitem{Krantz11}  S. Krantz, {\it A Mathematician Comes of Age}, Mathematical Association of America, 2011.

\bibitem{Mumford00} D. Mumford, The Dawning of the Age of Stochasticity, chapter in {\it Mathematics: Frontiers and Perspectives}, edited by V.Arnold, M.Atiyah, P.Lax and B.Mazur, AMS, 2000.

\bibitem{Ruelle02} D. Ruelle, Dynamical Zeta
Functions and Transfer
Operators, {\it Notices of the American Mathematical Society}, volume 49, issue 8 (2002), pp. 887-895.



\bibitem{Tao05}  T. Tao, The Dichotomy between Structure and Randomness, Arithmetic Progressions, and the Primes,  http://lanl.arxiv.org/pdf/math/0512114.pdf.


	







\end{thebibliography}
\end{document}